\documentclass{svproc}
\usepackage{url}

\usepackage{hyperref}
\usepackage{type1cm}        
\usepackage{makeidx}         
\usepackage{graphicx}        

\usepackage{multicol}        
\usepackage[bottom]{footmisc}

\usepackage{newtxtext}       %
\usepackage[varvw]{newtxmath}       

\usepackage{bm}
\usepackage{siunitx}
\usepackage{xcolor}
\usepackage{amsmath,bm}
\DeclareMathOperator{\E}{\mathbb{E}}
\usepackage{booktabs}
\usepackage{subcaption}
\usepackage{multirow}
\usepackage{varwidth}
\usepackage{tabularx}
\usepackage[export]{adjustbox}
\newcommand{\norm}[1]{\left\lVert#1\right\rVert}
\usepackage{lscape}
\newcommand*{\doi}[1]{\href{https://doi.org/\detokenize{#1}}{doi: \detokenize{#1}}}
\newcounter{algorithm}

\DeclareSymbolFont{bbold}{U}{bbold}{m}{n}
\DeclareSymbolFontAlphabet{\mathbbold}{bbold}

\usepackage[algo2e]{algorithm2e} 
\usepackage{caption,lipsum}
\usepackage{algorithm,algpseudocode}
\SetArgSty{textnormal}

\begin{document}

\mainmatter              
\title{Application of quasi-Monte Carlo in Mine Countermeasure Simulations with a Stochastic Optimal Control Framework}
\titlerunning{Application of quasi-Monte Carlo in Mine Countermeasure Simulations with a Stochastic Optimal Control Framework}  
%
\author{Philippe Blondeel\inst{1}, Filip Van Utterbeeck\inst{1} \and Ben Lauwens\inst{1} }
\authorrunning{Philippe Blondeel et al.} 
%
\tocauthor{Ivar Ekeland, Roger Temam, Jeffrey Dean, David Grove,
Craig Chambers, Kim B. Bruce, and Elisa Bertino}
\institute{Royal Military Academy,Avenue de la Renaissance 30, 1000 Brussels, Belgium\\
\email{\{philippe.blondeel, filip.vanutterbeeck, ben.lauwens\}@mil.be}
}

\titlerunning{Application of quasi-Monte Carlo in Mine Countermeasure Simulations}
\maketitle              

\begin{abstract}
Modelling and simulating mine countermeasures search missions performed by autonomous vehicles equipped with a sensor capable of detecting mines at sea is a challenging endeavour. The output of our stochastic optimal control implementation consists of an optimal trajectory in a square domain for the autonomous vehicle such that the total mission time is minimized for a given residual risk of not detecting sea mines. We model this risk as an expected value integral. We found that upon completion of the simulation, the user requested residual risk is usually not satisfied. We solved this by implementing a relaxation strategy which consists of incrementally increasing the square search domain. We then combined this strategy with different quasi-Monte Carlo schemes used for solving the integral. We found that using  a Rank-1 Lattice scheme yields a speedup up to a factor two with respect to the Monte Carlo scheme.  We also present an implementation which allows us to compute a trajectory in a convex quadrilateral domain, as opposed to a square domain, and combine it with our relaxation strategy.

\keywords{Stochastic Optimal Control, Mine Countermeasures, quasi-Monte Carlo}
\end{abstract}
\section{Introduction}
Modelling and simulating mine countermeasures (MCM) search missions performed by autonomous vehicles is a challenging endeavour. The goal of these simulations typically consists of computing an optimal trajectory for the autonomous vehicle in a designated zone where the presence of underwater sea mines is suspected. This trajectory is such that the residual MCM risk in the zone is below a certain threshold, while it is simultaneous ensured that the mission time is below a certain value.  This type of problem is referred to as a Coverage Path Planning (CPP) problem, see \cite{Abreu,AI2021110098,https://doi.org/10.1049/rsn2.12256,HowieChoset,GALCERAN20131258,1511771,6380733,9121686}.  From the previously cited works, two main approaches emerge. The first one consists of dividing  the to-be surveyed zone into a grid consisting of square or hexagonal cells. This grid is then used to compute a trajectory for the autonomous vehicle such that each cell is visited at least once. The second approach consists of computing a trajectory according to a pattern such as, a boustrophedon pattern, a zigzag search pattern or a weaving pattern. A novel third approach is presented in the work of  \cite{9140316}, and implemented for zones which consist of square search domains in our previous work \cite{eccomas_Blondeel} in a stochastic optimal control  framework \cite{pulsipher2022unifying}. An important constraint in our implementation consists of the MCM residual risk constraint. This constraint is formulated as an expected value integral over the square search domain. It has as purpose to ensure that the residual risk of mines not having been detected  is below a certain user-defined threshold after completion of the simulated mission in said domain. In order to numerically evaluate this integral, a Monte Carlo (MC) or a quasi-Monte Carlo (qMC) integration scheme is used. However, we noticed that upon completion of the simulation that the numerically calculated residual risk value is often higher than the requested one. In this work, we therefore propose a solution to this issue by implementing a relaxation strategy of the stochastic optimal control problem. Our strategy consists of successively increasing the considered square search domain until the user requested residual risk is satisfied.  We combine and compare our implementation with different qMC point sets such as a Rank-1 Lattice rule, a Sobol' and interlaced Sobol'  sequences. In addition to this, we also provide an implementation based on the generation of qMC points for triangular domains, which enables us to compute an optimal trajectory for the autonomous vehicle when a convex quadrilateral domain is considered. We then combine this implementation with our relaxation strategy.

The paper is structured as follows. First we present the methodology, where we introduce the sensor model, the formulation of the stochastic optimal control problem, elaborate upon how we numerically evaluate the residual risk integral, introduce our relaxation strategy, and discuss how to compute an optimal trajectory in a convex quadrilateral domain. Second, we present the results of our relaxation strategy combined with the aforementioned qMC points, and compare them against themselves and a MC point set. We also present results for the convex quadrilateral domain. Last, we offer a conclusion, and present directions for future research.

\section{METHODOLOGY}

In this section, we first present the equations describing the sensor model for a Forward Looking Sensor (FLS), i.e., a sensor which is only capable of detecting the presence of sea mines when the mines  are located in front of the autonomous vehicle. Hereafter, we give the equations describing the stochastic optimal control problem as set forward in \cite{9140316}. Third, we briefly discuss how we numerically evaluate the residual risk integral. Fourth, we present our relaxation strategy and provide an algorithm. Last, we present our implementation which enables us to compute trajectories in a convex quadrilateral domain as opposed to square domains.

\subsection{The sensor model}

We briefly introduce the equations pertaining to the modelisation of a FLS. A more detailed description can be found in \cite{9140316}.

The sensor model is given by
\begin{equation}
\gamma\left(\bm{x}\left(t\right),\bm{\omega}\right) := \lambda\,p\left(\bm{x}\left(t\right),\bm{\omega}\right)\,F_\alpha\left(\bm{x}\left(t\right),\bm{\omega}\right)\,F_\varepsilon\left(\bm{x}\left(t\right),\bm{\omega}\right),
\label{eq:sensor_model}
\end{equation}
where $\lambda$ stands for the Poisson scan rate in ${\text{s}}^{-1}$.

The sensor model of Eq.\,\eqref{eq:sensor_model} consists of three parts. The first part consists of $p\left(\bm{x}\left(t\right),\bm{\omega}\right)$ and is given by
\begin{equation}
p\left(\bm{x}\left(t\right),\bm{\omega}\right) := \Phi_\text{n}\left(\frac{\text{FOM} - 20 \log_{10}\left(\norm{\bm{\omega} - \bm{x}\left(t\right)} + a\norm{\bm{\omega}-\bm{x}\left(t\right)}\right)}{\sigma}\right),
\end{equation}
where  $\Phi_\text{n}\left(\cdot\right)$ stands for the cumulative density function (CDF) of the normal distribution, $\text{FOM}$ is a parameter related to the sonar characteristics, $a$ is the attenuation coefficient in dB/km, and $\bm{x}\left(t\right)$ is given in Eq.\,\eqref{eq:pos}.
We give $\norm{\bm{\omega} - \bm{x}\left(t\right)}$ as
\begin{equation}
\norm{\bm{\omega} - \bm{x}\left(t\right)} := \sqrt{\left(\omega_x - x\left(t\right)\right)^2} + \sqrt{\left(\omega_y - y\left(t\right)\right)^2},
\end{equation}
where $\bm{\omega}$ is defined as  $\bm{\omega}:= \left(\omega_x,\omega_y\right)$, and stands for the position of the possible target, i.e., a sea mine.
This second part, $F_\alpha\left(\bm{x}\left(t\right),\bm{\omega}\right)$, models the detection of the sensor in front of  the autonomous vehicle according to its Field Of View (FOV), and is given as
\begin{equation}
F_\alpha\left(\bm{x}\left(t\right),\bm{\omega}\right):=\frac{1}{1+e^{p_\alpha\left(-\frac{\alpha_{\text{FOV}}}{2} - \alpha^b\left(\bm{x}\left(t\right),\bm{\omega}\right)\right)}} + \frac{1}{e^{p_\alpha\left(\alpha^b\left(\bm{x}\left(t\right),\bm{\omega}\right)-\frac{\alpha_{\text{FOV}}}{2}\right)}} - 1,
\label{eq:F_a}
\end{equation}
where $p_\alpha$ is a parameter used to adjust the slope of the sigmoidal curves, and
\begin{equation}
  \begin{gathered}
  \alpha^b\left(\bm{x}\left(t\right),\bm{\omega}\right) := \arctan_2\left(dx^b\left(\bm{x}\left(t\right),\bm{\omega}\right),dy^b\left(\bm{x}\left(t\right),\bm{\omega}\right)\right) \\
dx^b\left(\bm{x}\left(t\right),\bm{\omega}\right):=\left(\omega_x -x(t)\right)\,\cos(\psi(t))+\left(\omega_y-y(t)\right)\,\sin(\psi(t)) \\
dy^b\left(\bm{x}\left(t\right),\bm{\omega}\right):=-\left(\omega_x -x(t)\right)\,\sin(\psi(t))+\left(\omega_y-y(t)\right)\,\cos(\psi(t)).
\end{gathered}
\label{eq:F_a_expanded}
\end{equation}

%
%
 In Eq.\,\eqref{eq:F_a},  $\alpha_\text{FOV}$ is the Field Of View angle of the FLS in degrees. In Eq.\,\eqref{eq:F_a_expanded}, $\psi(t)$ represents the angle the autonomous vehicle has with respect to the horizontal axis in degrees, $x(t)$ is the position along the $x$-axis and $y(t)$ is the position along the $y$-axis.
 
The third part, $F_\varepsilon\left(\bm{x}\left(t\right),\bm{\omega}\right)$ accounts for the height $h$ in meters above the ocean floor of the sensor attached to the autonomous vehicle, and is  given by
\begin{equation}
F_\varepsilon\left(\bm{x}\left(t\right),\bm{\omega}\right):=\frac{1}{1+e^{p_\epsilon\left(\varepsilon_\text{DE}-\frac{\varepsilon_\text{FOV}}{2} - \varepsilon^b\left(\bm{x}\left(t\right),\bm{\omega}\right)\right)}} + \frac{1}{e^{p_\epsilon\left(\varepsilon^b\left(\bm{x}\left(t\right),\bm{\omega}\right)-\varepsilon_\text{DE}-\frac{\varepsilon_\text{FOV}}{2}\right)}} - 1,
\end{equation}
where $p_\varepsilon$ is a parameter used to adjust the slope of the sigmoidal curves, and
\begin{equation}
\varepsilon^b\left(\bm{x}\left(t\right),\bm{\omega}\right) := \arctan\left(\frac{-h}{\norm{\bm{\omega} - \bm{x}\left(t\right)} }\right),
\end{equation}
where $\varepsilon_\text{FOV}$ and $\varepsilon_\text{DE}$ respectively stand for the vertical FOV, and the downward elevation angle such that the sensor can ensonify the sea floor. For a more thorough description we refer to \cite{9140316}.

\subsection{Formulation of the stochastic optimal control problem}
\label{sec:Trajectories}

In \cite{9140316}, the optimization problem was formulated such that the residual MCM risk is minimized for a given mission time. The residual MCM risk is defined as the probability that one or more autonomous vehicles fail to detect the mines in a search domain by the end of an MCM operation. The residual MCM risk is sometimes also referred to as `the probability of non-detection'.

 Here and in our previous work \cite{eccomas_Blondeel},  the optimization problem is formulated such that the  mission time $T_f$ needed to survey a designated domain $\Omega$ is minimized for a given residual MCM risk,
\begin{equation}
\text{min}\, T_f,
\end{equation}
subjected to
\begin{equation}
 \E[q\left(T_F\right)] :=  \int_\Omega \text{e}^{-\int_0^{T_F} \gamma\left(\bm{x}\left(\tau\right),\bm{\omega}\right)\, d\,\tau}\phi\left(\bm{\omega}\right) d\,\bm{\omega} \leq  \text{residual MCM risk}
\label{eq:exp}
\end{equation}
where $\phi(\cdot)$ is the probability density function (PDF) of the  distribution against which we integrate.  In this case we consider $\phi(\cdot)$ to be the PDF of the uniform distribution. This is because we consider a uniform distribution of the possible targets, i.e., the sea mines. In \S\,\ref{sec:ResdiualRisk}, we will give a short overview of how this integral is evaluated numerically.

The position of the autonomous vehicle is given by 

\begin{equation}
\bm{x}\left(t\right) := f(x(t), y(t), \psi(t), r(t)),
\label{eq:pos}
\end{equation}
with $r(t)$ being the turn rate in degrees per second.
The position at time $t$, see Eq.\,\eqref{eq:pos}, is governed by the following differential equations,
\begin{equation}
\begin{aligned}
&\frac{d\,x(t)}{d\,t} = V \text{cos}(\psi(t)) \\
&\frac{d\,y(t)}{d\,t} = V \text{sin}(\psi(t)) \\
&\frac{d\,\psi(t)}{d\,t} = r(t) \\
&\frac{d\,r(t)}{d\,t} = \frac{1}{T}\left(K p(t) - r(t)\right),
\end{aligned}
\end{equation}
where $V$ stands for the speed in m/s, $K$ is the Nomoto gain constant with units $\text{s}^{-1}$, $T$ is the Nomoto time constant with units s, and $p(t)$ is the rudder deflection angle given in degrees. The optimization software \cite{pulsipher2022unifying} will control the values for the rudder deflection angle in order to compute a trajectory satisfying the target function and the constraints. The parameters which will be used for all simulations are listed in Tab.\,\ref{Tab:Param}.

\begin{table}
\hspace{-1.cm}
{\renewcommand{\arraystretch}{1.2}
\setlength{\tabcolsep}{4pt}
\scalebox{0.85}{
\begin{tabular}{cccccccccccccc}
\toprule
 {Parameter name} &  $\alpha_{\text{FOV}}$ &$h$ & $\sigma$ & $\lambda$ &$a$ & $\varepsilon_{\text{FOV}}$ &$\varepsilon_{\text{DE}}$ & $\text{FOM}$ &$p_\alpha$ &$p_\varepsilon$ & $\text{V}$ &$\text{T}$ & $\text{K}$  \\
  \cmidrule(rl{4pt}){1-14}  
 {Value} & $120.0^\text{o}$ &20.0\,m &9.0\,[/]&20.0\,$\text{s}^{-1}$&5.2\,dB/km&$5.0^\text{o}$&$-6.0^\text{o}$&72.0\,[/]&25.0\,[/]&400.0\,[/]&2.5\,m/s&0.5\,s&5.0\,$\text{s}^{-1}$\\
\bottomrule
\end{tabular}}}
\caption{Parameters for the simulations.}
\label{Tab:Param}
\end{table} 
\raggedbottom

\begin{figure}
\centering
\includegraphics[scale=0.22]{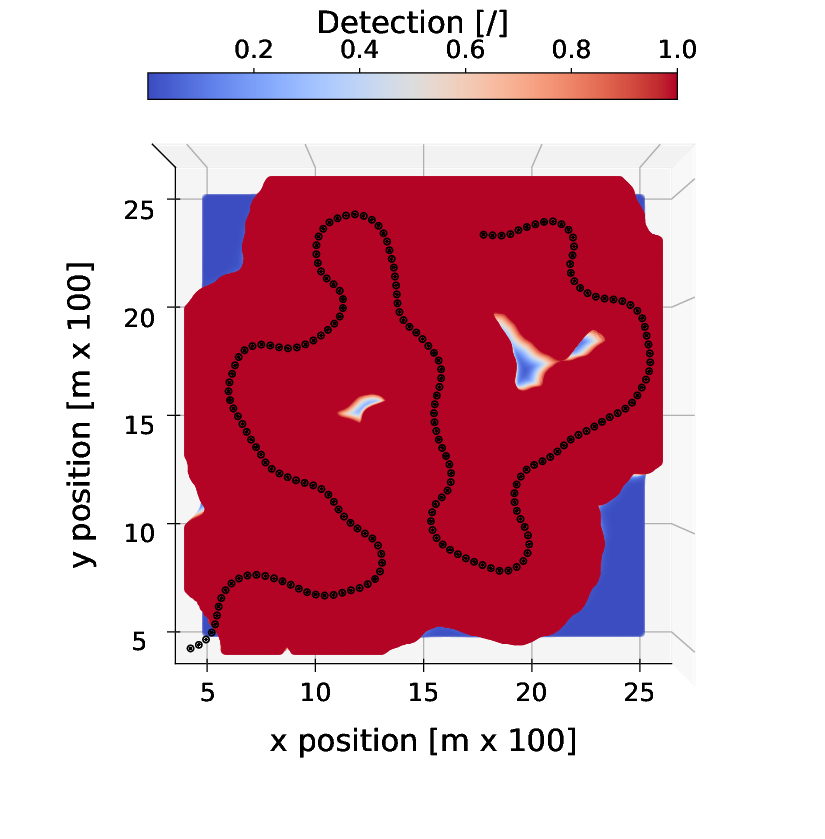}
\vspace{-0.5cm}
\caption{A typical solution produced by our stochastic optimal control implementation of the MCM problem when considering an autonomous vehicle in the $\left[5,25\right]^2$ domain with a requested residual risk of maximally $10\,\%$. The black dotted line represent the trajectory, the red area is the part of the domain that has been `seen' by the onboard sensor, and the blue area is the part that has not been `seen' by the sensor.}
\label{fig:example_sol}
\end{figure}

In Fig.\,\ref{fig:example_sol}, we show an example of a typical  solution of the stochastic optimal control problem. We consider a square domain $\Omega \in [5,25]^2$ with a requested residual risk of maximally $10\,\%$. The optimal trajectory is represented by the black dotted lines. The area of the domain which is coloured red has been `seen' (ensonified) by the sensor attached to the autonomous vehicle. The area in blue has not been ensonified by the sensor while following the depicted trajectory. 

\subsection{Residual MCM risk and numerical implementation of the problem}
\label{sec:ResdiualRisk}
In our implementation, the residual MCM risk integral as defined in Eq.\,\eqref{eq:exp} is approximated by two possible  sampling schemes. The first approach consists of a MC scheme where we compute the integral as follows,
\begin{equation}
\int_\Omega \text{e}^{-\int_0^{T_F} \gamma\left(\bm{x}\left(\tau\right), \bm{\omega}\right)\, d\,\tau}\phi\left(\bm{\omega}\right) d\,\bm{\omega} \approx \frac{1}{N}\sum_{n=1}^N \text{e}^{-\int_0^{T_F} \gamma\left(\bm{x}\left(\tau\right),     {\text{\textbf{x}}}^{(n)}\right)\, d\,\tau},
\end{equation}

where ${\text{\textbf{x}}}^{(n)}$ are MC samples originating from the uniform distribution such that ${\text{\textbf{x}}}^{(n)} \in \Omega$, and $N$ stands for the number of samples.
The second  approach consists of a qMC integration scheme where the qMC points are distributed uniformly in the unit cube. For integration over  $\mathbb{R}^s$  against the uniform distribution, the most commonly used approach, see \cite{Graham,KUO}, consists of performing a change of variables on the multidimensional integral, i.e., $\boldsymbol{y} = \Phi^{-1}(\boldsymbol{\zeta})$, such that the integration domain is changed from $\mathbb{R}^s$ to $\left[0,1\right]^s$,
 \begin{equation}
 \int_{\mathbb{R}^s} f\left(\boldsymbol{y}\right) \phi(\boldsymbol{y})\, \text{d}\boldsymbol{y} = \int_{[0,1]^s} f(\Phi^{-1}(\boldsymbol{\zeta}) ) \text{d} \,\boldsymbol{\zeta},
 \label{eq:inv_g}
 \end{equation}
 where ${\Phi}^{-1}$ is the inverse univariate uniform cumulative distribution function.  The right-hand side of Eq.\,\eqref{eq:inv_g} can then be evaluated numerically with the shifted qMC quadrature points $\mathbf{u}^{\left(r,n\right)}$  defined in the unit cube. Applied to the residual risk integral, the qMC formulation is as follows, 
\begin{equation}
\int_\Omega \text{e}^{-\int_0^{T_F} \gamma\left(\bm{x}\left(\tau\right), \bm{\omega}\right)\, d\,\tau}\phi\left(\bm{\omega}\right) d\,\bm{\omega} \approx \frac{1}{R}\sum_{r=1}^{R}\frac{1}{N}\sum_{n=1}^N \text{e}^{-\int_0^{T_F} \gamma\left(\bm{x}\left(\tau\right), \Phi^{-1}\left(\text{\textbf{u}}^{(n,r)}\right)\right)\, d\,\tau},
\label{eq:qmc_impl}
\end{equation}
where the inverse of the univariate uniform cumulative distribution function is applied  point-wise to the qMC points, and $R$ stands for the number of shifts.

Before describing the issue arising with our initial `naive' implementation of the mine counter measure simulation, we first give a short description of  how we assess the residual MCM risk at completion of the simulation. For a given domain $\Omega$ and starting coordinate for the survey vehicle $\boldsymbol{\xi}$, we request a maximal allowable residual MCM risk. This data is then fed into the stochastic optimal control implementation, which internally evaluates the residual MCM risk integral with our supplied MC or qMC points. Upon completion of the simulation,  the position of the autonomous vehicle at discretized time points, see Eq.\,\eqref{eq:pos} is returned. With this data, we recompute, in a post-processing step, the MCM risk integral as defined in Eq.\,\eqref{eq:exp} with a Monte Carlo integration scheme which uses $2^{14}$ points on the $\Omega$ domain. The value resulting from this post-processing computation is then used to assess if the requested residual MCM risk is satisfied or not. An  overview of the described implementation is given in Algorithm\,\ref{algo:first_implement}.

The arising issue is such that at the end of our simulation, the MCM risk computed in post-processing was often higher than the maximally allowed user defined risk. We illustrate this issue by solving  $1000$ independent simulation runs of our problem with qMC and MC. We show the individual results, and the boxplots in Fig.\,\ref{fig:residual_mcm_risk} and Tab.\,\ref{Tab:residual_mcm_no_increase}. We observe for  both schemes,  that the user requested residual MCM risk ($5\%$) is typically not satisfied at the end of the optimisation process. When planning mine counter measure mission, this erroneous result might lead to a false sense of security. Furthermore, we also observe that the implementation with a qMC scheme yields better results than the MC one. The mean and median are lower for the qMC case, and the results are less `spread out'. This is reflected in a lower standard deviation for the qMC case.

For the generation of qMC points we used a Rank-1 Lattice rule  constructed from a 10-dimensional base-2 generating vector for an equal weight Korobov space from \cite{Hickernell2012}. The text file can be found on \cite{HickernellVec}.

\begin{table}[ht]
\begin{minipage}[b]{0.5\linewidth}
\hspace{-1.8cm}
\includegraphics[scale=0.65]{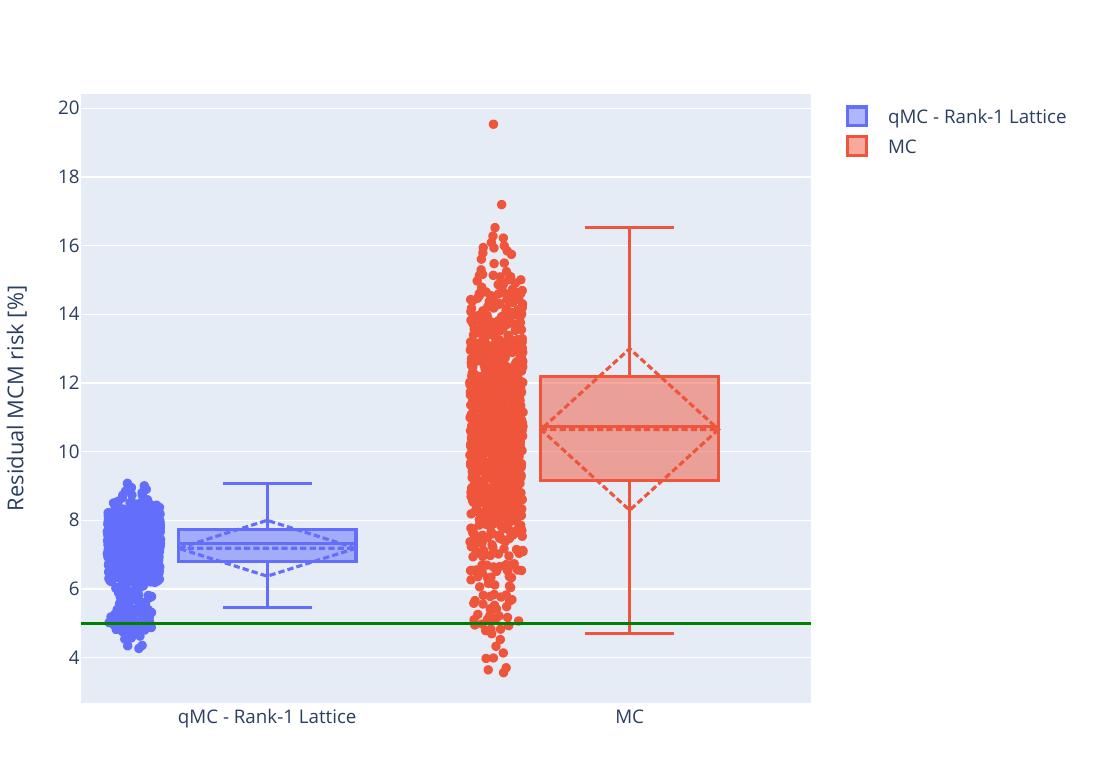}
\captionof{figure}{Boxplot and individual results of the residual MCM risk in \% of 1000 simulations. The green line indicates the (maximal allowable) 5\% user requested residual MCM risk which all individual simulations should satisfy.}
\label{fig:residual_mcm_risk}
\end{minipage}
\hfill
\begin{minipage}[b]{0.4\linewidth}
\vspace{-6cm}
\scalebox{1.0}{
\begin{tabular}{ccc}
\toprule
 {} &  {qMC } &   \multirow{2}*{MC} \\
 {} &  {Rank-1 Lattice} &    \\
 \cmidrule(rl{4pt}){1-3}  
 Median & 7.33\,\% &   10.73\,\%   \\
 Mean & 7.18\,\% & 10.64\,\%    \\
 Standard deviation & 0.81\,\%&  2.35\,\%        \\
 Minimum &  4.27\,\%  & 3.56\,\%\\
 Maximum &  9.07\,\% & 19.53\,\% \\
\bottomrule
\end{tabular}}
\vspace{3.cm}
\caption{Numerical results belonging to Fig.\,\ref{fig:residual_mcm_risk}}
\label{Tab:residual_mcm_no_increase}
  \end{minipage}%
\end{table} 

\begin{algorithm2e}
 \caption{Implementation of the stochastic control problem.}
 \KwData{Domain $\Omega$, Starting points for survey vehicle $\boldsymbol{\xi}$, request MCM risk $M$, MC or qMC point set $\text{\textbf{k}}$}
Solve the stochastic optimal control problem as defined in \S\,\ref{sec:Trajectories} on $\Omega$\ with $M$, $\text{\textbf{k}}$, and $\boldsymbol{\xi}$\;
Compute $\int_\Omega \text{e}^{-\int_0^{T_F} \gamma\left(\bm{x}\left(\tau\right),\bm{\omega}\right)\, d\,\tau}\phi\left(\bm{\omega}\right) d\,\bm{\omega}$ with $2^{14}$ MC points \;
$\text{MCM Risk} \gets \int_\Omega \text{e}^{-\int_0^{T_F} \gamma\left(\bm{x}\left(\tau\right),\bm{\omega}\right)\, d\,\tau}\phi\left(\bm{\omega}\right) d\,\bm{\omega}$ \;
\vspace{0.1cm}
\Return MCM Risk
\vspace{0.4cm}
 \label{algo:first_implement}
\end{algorithm2e}

\subsection{Relaxation strategy for the residual MCM risk}
In order to tackle  the issue  put forward in \S\,\ref{sec:ResdiualRisk}, we developed a relaxation strategy under the form of an adaptive domain increase algorithm. This algorithm ensures that the  MCM risk computed in a post-processing step is lower or at most equal to the requested residual MCM risk. We give the implementation in Algorithm\,\ref{algo:adaptive}. The main idea consists of sequentially solving the stochastic optimal control problem with an ever increasing size of the domain. The qMC or MC points we use are generated once at the start of the simulation on the $\left[0,1\right)^2$ domain, after which they are mapped to the desired domain. Important to note is that when we combine our relaxation strategy with MC points, we use the same mapping technique as the one for the qMC points, see \S\,\ref{sec:ResdiualRisk}. When using qMC points, we use one shift, i.e., $R = 1$ in Eq.\,\eqref{eq:qmc_impl}. We observed that using more than one shift did not yield a further advantage. The total computational cost in seconds is computed as the sum of the computational costs needed for each solve of the stochastic optimal control problem.  The values for the two dimensional vector $\Xi$, in Algorithm\,\ref{algo:adaptive}, govern the domain increasing following $\Omega_{\text{new}} = \Omega_{\text{old}} + \Xi$. The values for $\Xi$ are taken small with respect to the size of the domain. In this case they are equal to 0.2 for both dimensions. We also note that the computation of the integral in the post-processing step is done with MC points generated on the initial domain ${\Omega_{\text{init}}}$.

\begin{algorithm2e}
 \caption{Adaptive implementation of the stochastic control problem.}
 \KwData{Domain $\Omega$, Starting points for survey vehicle $\boldsymbol{\xi}$, request MCM risk $M$, two-dimensional vector $\Xi$, MC or qMC point set $\text{\textbf{k}}$}
$\Omega_{\text{init}} \gets \Omega$ \;
\While{$ \textup{MCM Risk} > M$}{
Solve the stochastic optimal control problem as defined in \S\,\ref{sec:Trajectories} on $\Omega$\ with $M$, $\text{\textbf{k}}$, and $\boldsymbol{\xi}$;

Compute $\int_{\Omega_{\text{init}}} \text{e}^{-\int_0^{T_F} \gamma\left(\bm{x}\left(\tau\right),\bm{\omega}\right)\, d\,\tau}\phi\left(\bm{\omega}\right) d\,\bm{\omega}$ with $2^{14}$ MC points \;
$\text{MCM Risk} \gets \int_{\Omega_{\text{init}}} \text{e}^{-\int_0^{T_F} \gamma\left(\bm{x}\left(\tau\right),\bm{\omega}\right)\, d\,\tau}\phi\left(\bm{\omega}\right) d\,\bm{\omega}$ \;

\If{$\textup{MCM Risk} > M$}{%
  $\Omega \gets \Omega + \Xi $ \tcc{Increase the domain $\Omega$}\
}
}
\vspace{0.1cm}
\Return MCM Risk
\vspace{0.4cm}
 \label{algo:adaptive}
\end{algorithm2e}

\subsection{Quadrilateral convex domains}
\begin{figure}
\begin{minipage}{0.33\textwidth}
\hspace{-0.5cm}
\includegraphics[scale=0.38]{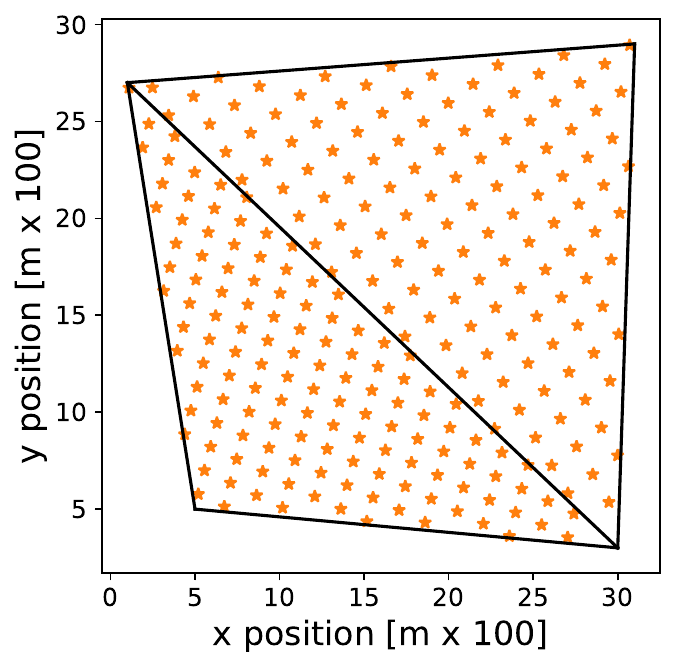}
\end{minipage}
\begin{minipage}{0.33\textwidth}
\hspace{-0.5cm}
\includegraphics[scale=0.38]{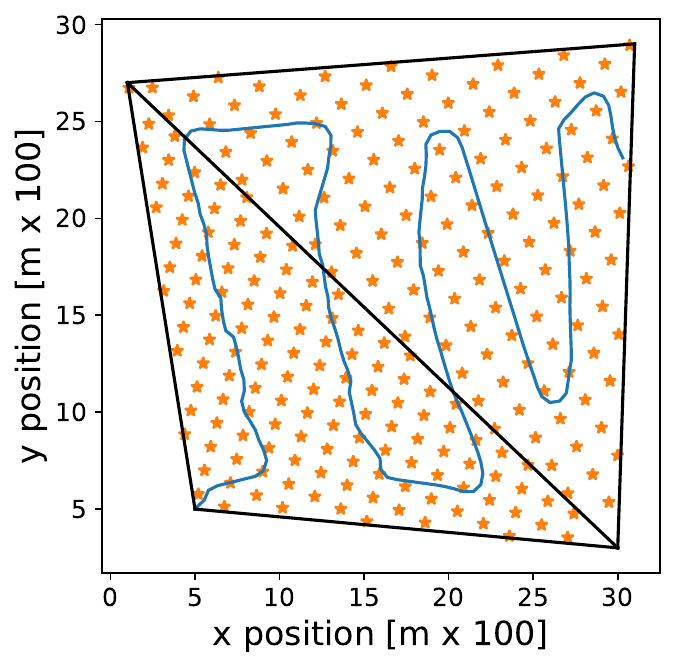}
\end{minipage}
\begin{minipage}{0.33\textwidth}
\hspace{-0.5cm}
\includegraphics[scale=0.20]{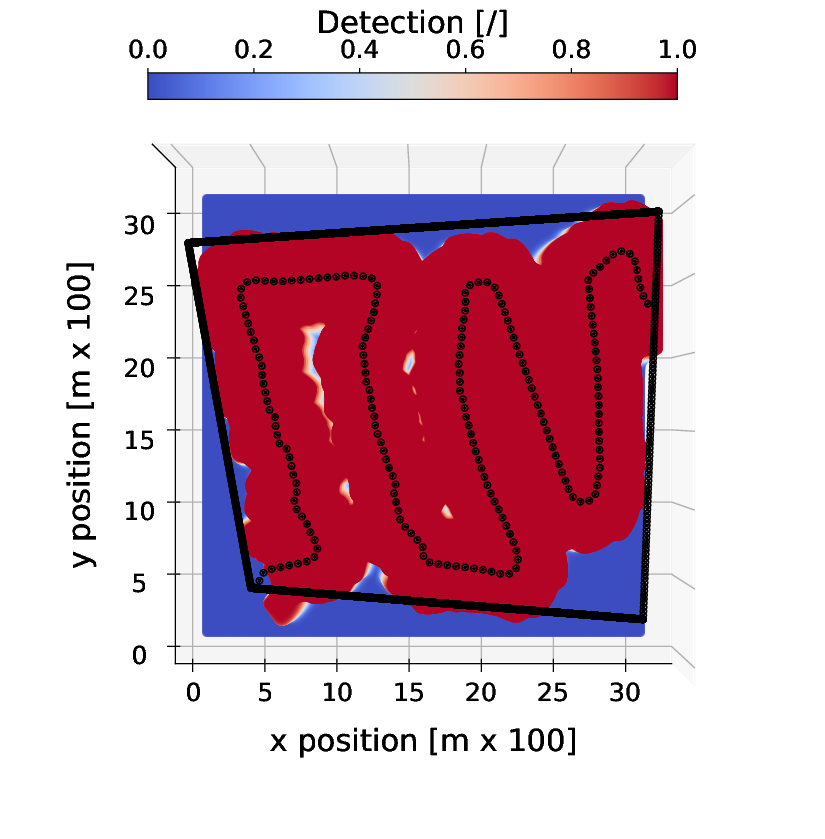}
\end{minipage}
\vspace{-0.8cm}
\caption{Using qMC points to compute a trajectory in a quadrilateral convex domain: (Left) Division of the convex quadrilateral domain into two triangles with the qMC points generated on each triangle in orange. (Center) The partitioned domain with the generated qMC points and the computed trajectory in blue. (Right) The domain delimited by the black lines, the ensonification of the domain in red and the trajectory represented by the black lines.}
\label{fig:quad_conv_dom}
\end{figure}

In addition to our relaxation strategy, we also developed an approach to compute trajectories in quadrilateral convex domains while using qMC points for the computation of the residual MCM risk integral of Eq.\,\eqref{eq:exp}. Our approach is as follows. First, we divide the quadrilateral convex domain into two distinct triangles. On each of these triangular domains we then generate qMC points for triangular domains according to the algorithm described in \cite{KinjalBasu}. The procedure is illustrated in Fig.\,\ref{fig:quad_conv_dom} where we show the partitioning of the convex quadrilateral domain, the generated qMC points on the triangles which make up the domain, and the ensonification of the domain.

\section{Results}
In this section we first present results when our relaxation strategy is applied to a square domain and combined with different qMC point sets and a MC point set. We consider  qMC point sets obtained from a Rank-1 Lattice rule, a Sobol' sequence, a Sobol' sequence with a factor two interlacing, and a Sobol' sequence with a factor three interlacing. In order to generate these qMC points, we use the \textbf{LatticeRules.jl} \cite{PieterJanGit3}  and \textbf{DigitalNets.jl} \cite{PhilippeBlondeelGit1} Julia implementations.  Second, we show results where we combine our relaxation strategy with a qMC point set for computing optimal trajectories in a convex quadrilateral domain. 
We show results pertaining to the residual risk computed in the post-processing step, the computation time, i.e., the time it takes to compute the solution, and the path time, i.e., the time needed for the autonomous vehicle to complete the survey of the domain. All our simulations have been performed on a workstation with an AMD 7542 processor containing 32 cores clocked at 2.96 GHz and 258 GB RAM.

\subsection{Relaxation strategy for square domains}
\begin{figure}
 \hspace{-1.cm}
\includegraphics[scale=0.35]{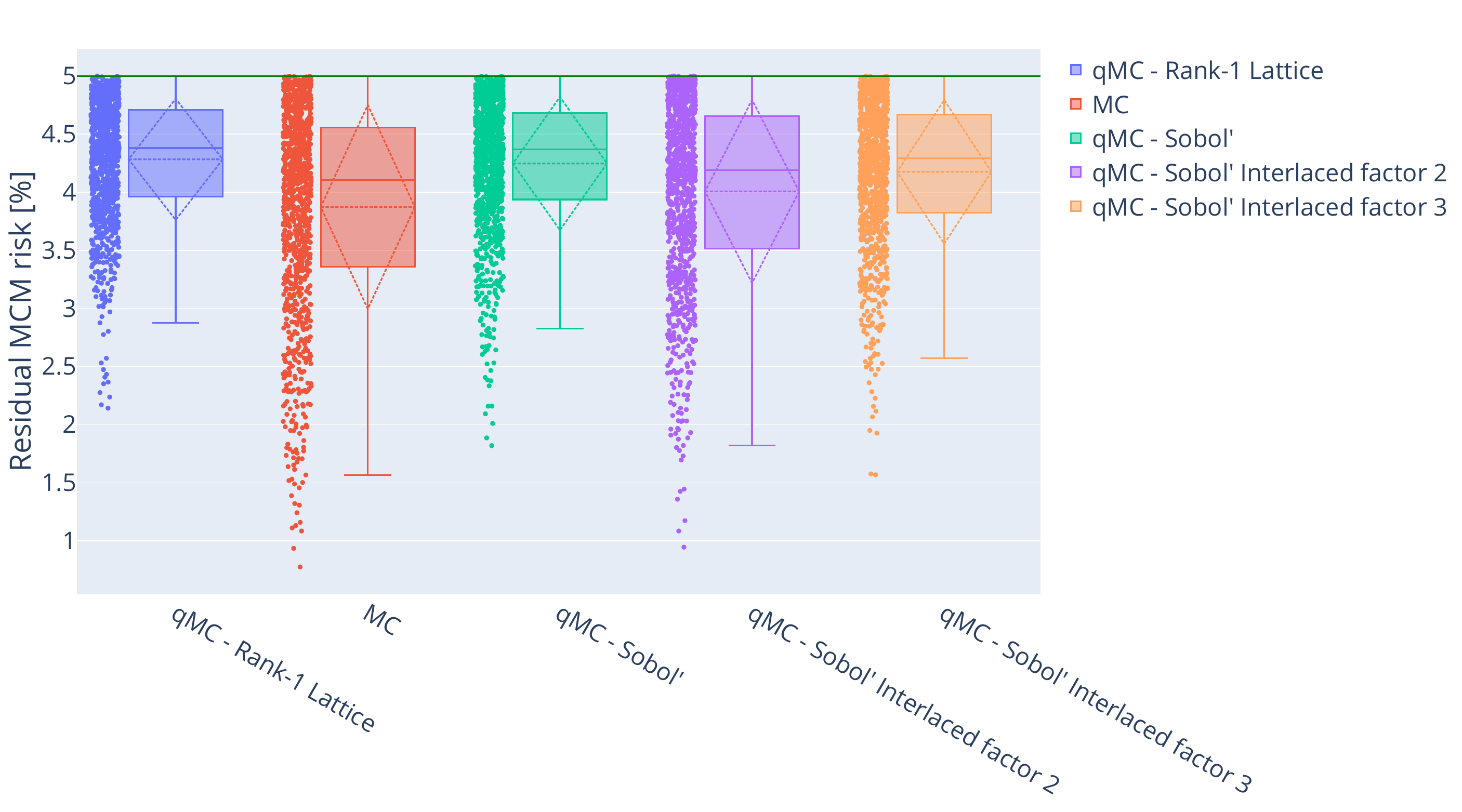}
\caption{The boxplots and individual results for the residual risk of 1000 independent simulations for different qMC and MC point sets. The green line at $5\,\%$ is the maximal user requested tolerance that all simulations should satisfy.}
\label{fig:residual_mcm_risk_adaptive}
\end{figure}

\begin{figure}
 \hspace{-1.cm}
\includegraphics[scale=0.35]{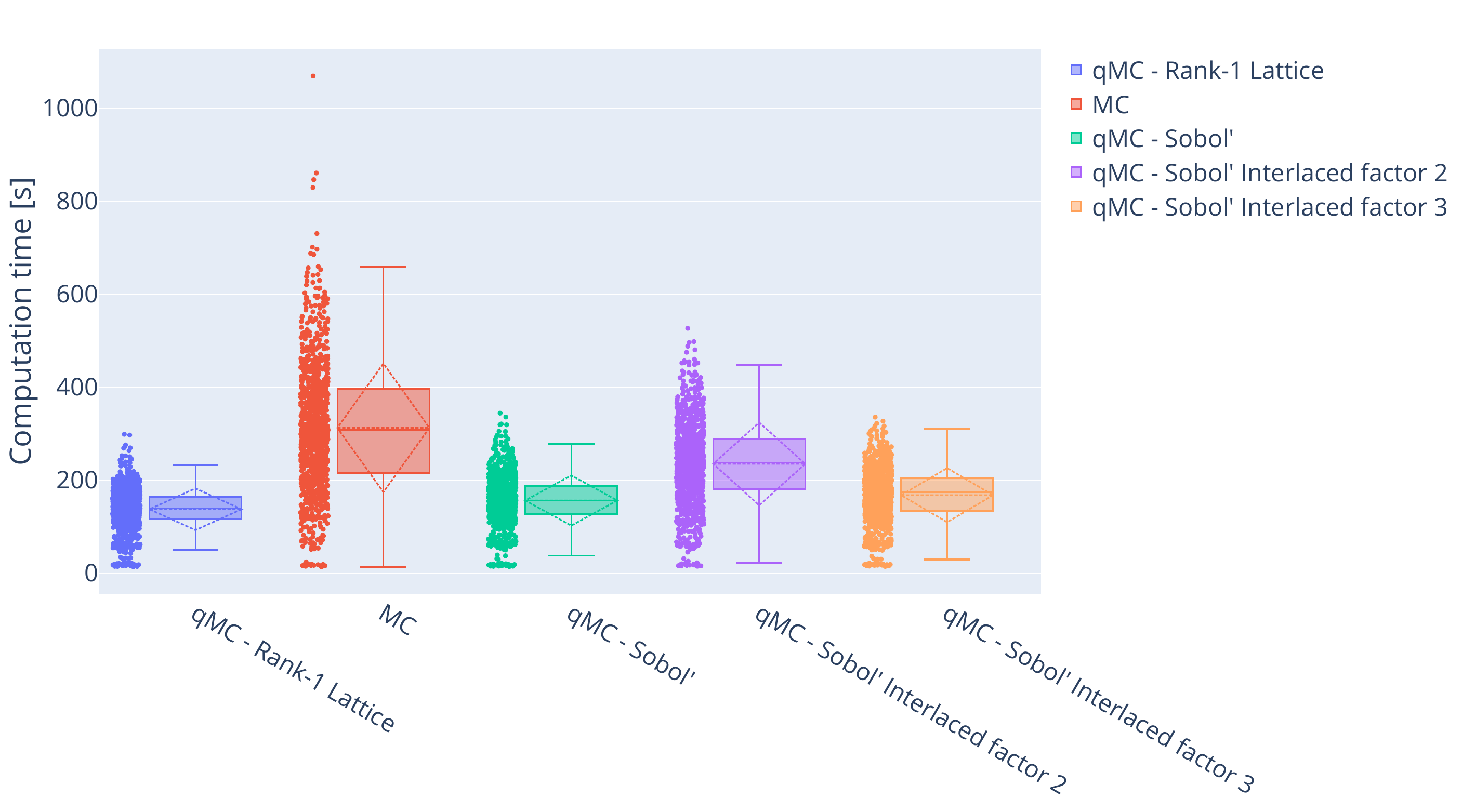}
\caption{The boxplots and individual results for the computation time of 1000 independent simulations for different qMC and MC point sets.}
\label{fig:adaptive_computation_time}
\end{figure}

\begin{figure}
 \hspace{-1.cm}
\includegraphics[scale=0.35]{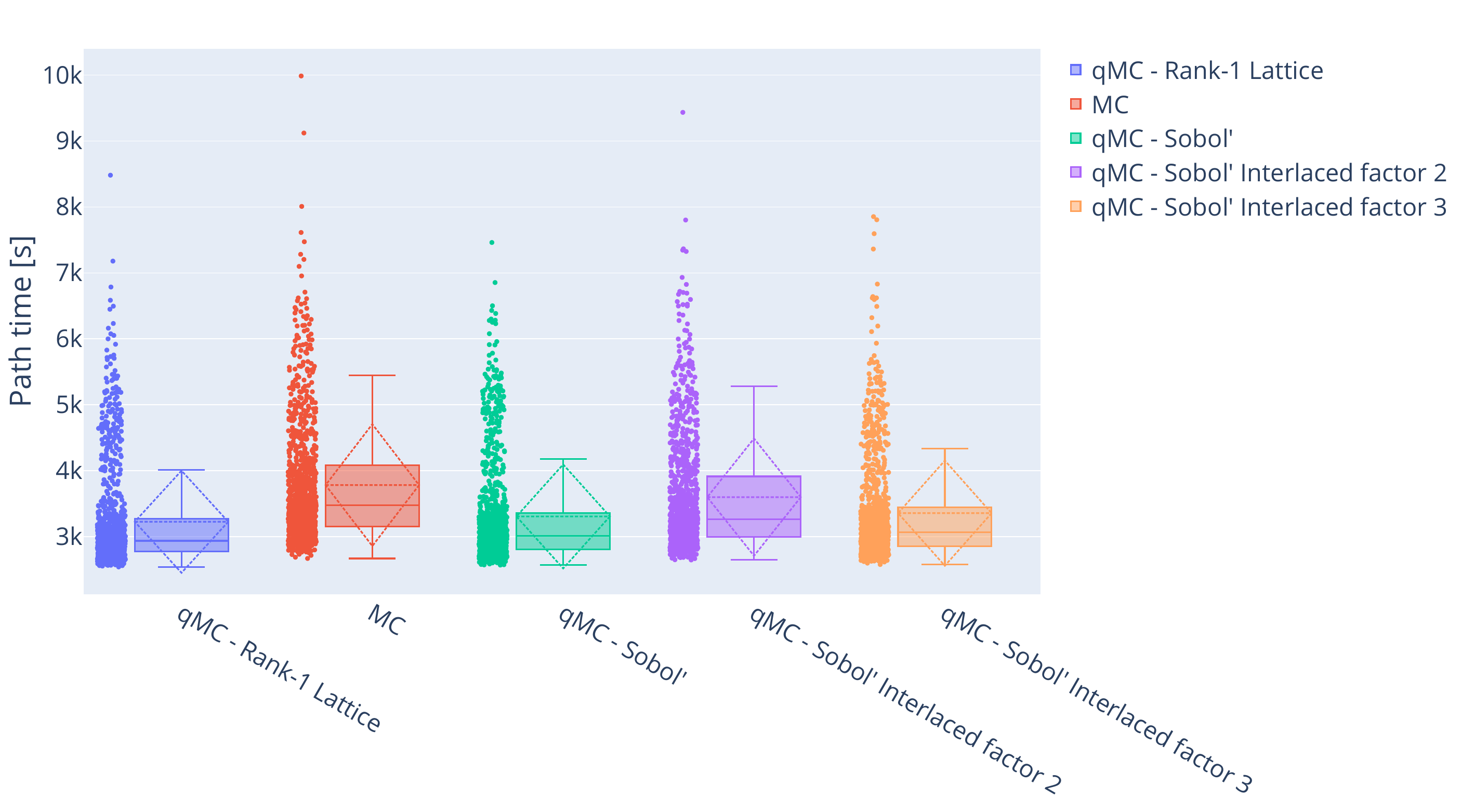}
\caption{The boxplots and individual results for the path time of 1000 independent simulations for different qMC and MC point sets.}
\label{fig:adaptive_run_time}
\end{figure}

In Fig.\,\ref{fig:residual_mcm_risk_adaptive} we present the stochastic results for 1000 independent simulations when applying our relaxation strategy to the stochastic optimal control problem. The median, mean, standard deviation, the minimum, and maximum are presented in Tab.\,\ref{Tab:residual_mcm_risk_adaptive}. We observe than when applying our relaxation strategy, the requested user tolerance is satisfied for all individual simulations. The computation times are presented in Fig.\,\ref{fig:adaptive_computation_time} and Tab.\,\ref{Tab:adaptive_computation_time}. From these results, we observe that when considering qMC point sets, the computation times are lower than when considering a MC point set. This is as expected based on our observations presented in \S\,\ref{sec:ResdiualRisk}. From all the considered qMC point sets, the Rank-1 Lattice rule is the one with the lowest associated computational cost.  We observe a speedup up to factor two with respect to MC.
In Fig.\,\ref{fig:adaptive_run_time} and Tab.\,\ref{Tab:adaptive_run_time} we present the path times. The shortest path times that satisfy the requested residual risks are associated with the qMC point sets and in particular with the Rank-1 Lattice rule. 

\begin{figure}
\hspace{-1.3cm}
\begin{minipage}{0.33\textwidth}
\includegraphics[scale=0.42]{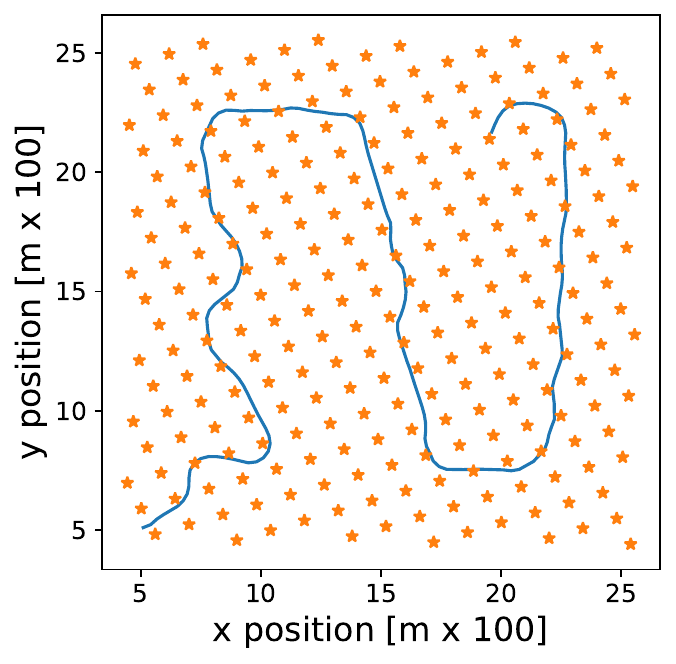}
\end{minipage}
\hspace{.5cm}
\begin{minipage}{0.33\textwidth}
\includegraphics[scale=0.42]{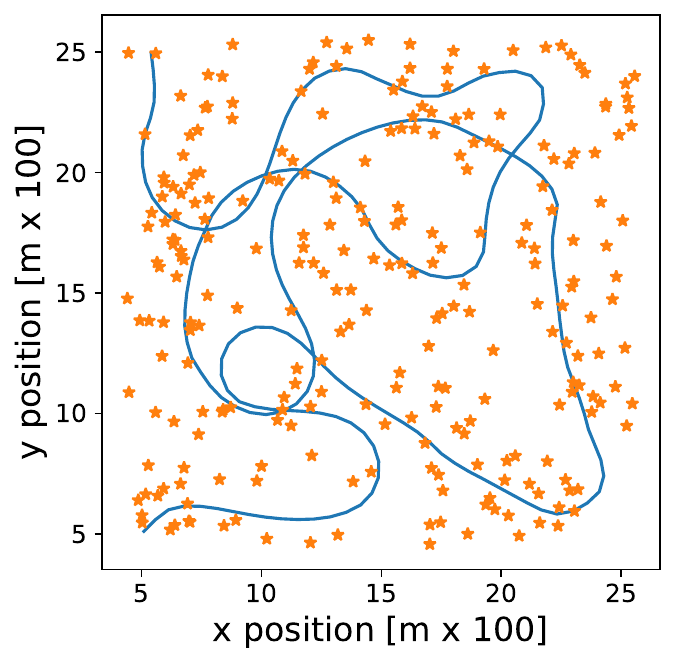}
\end{minipage}
\hspace{.5cm}
\begin{minipage}{0.33\textwidth}
\includegraphics[scale=0.42]{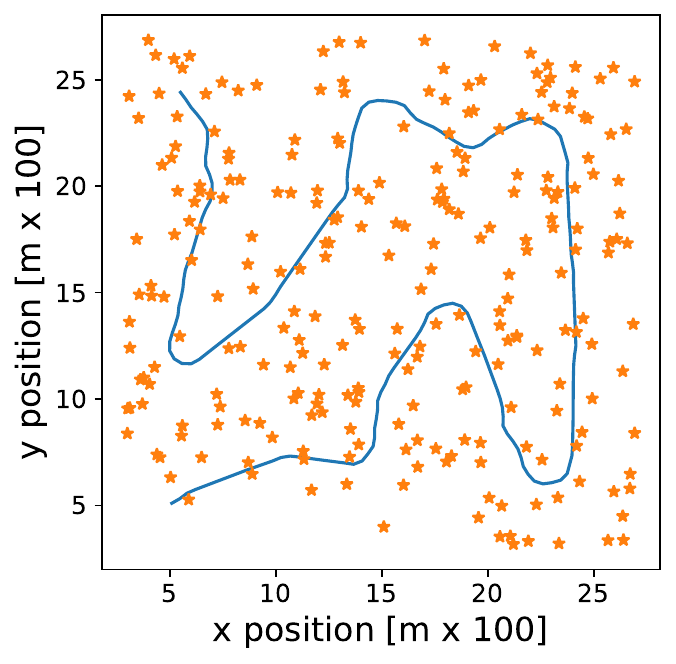}
\end{minipage}

\hspace{-2.0cm}
\begin{minipage}{0.33\textwidth}
\includegraphics[scale=0.22]{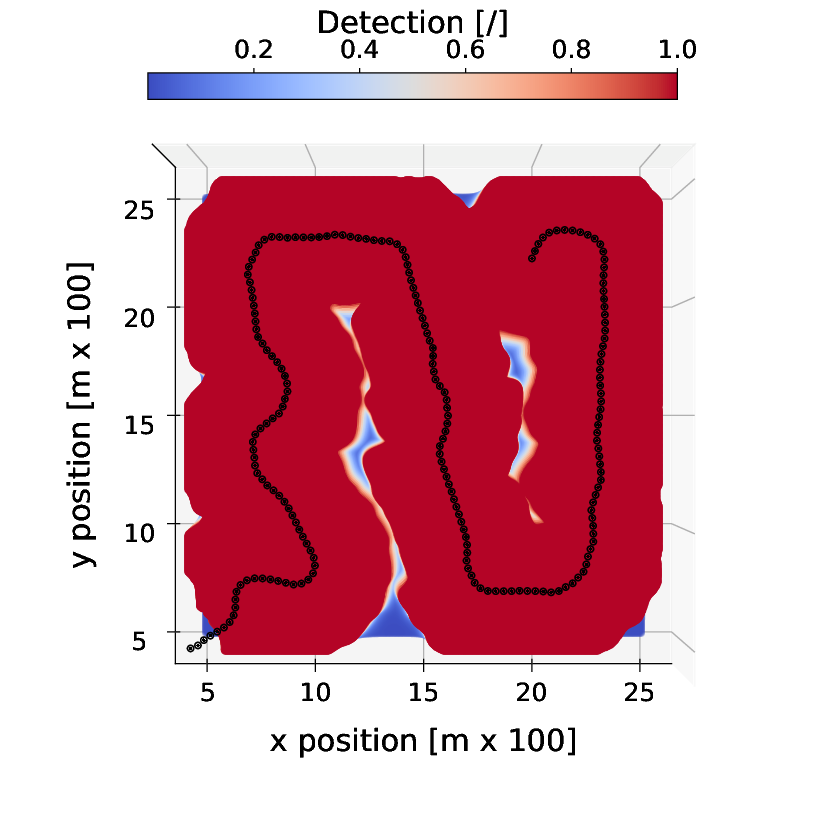}
\end{minipage}
\hspace{.5cm}
\begin{minipage}{0.33\textwidth}
\includegraphics[scale=0.22]{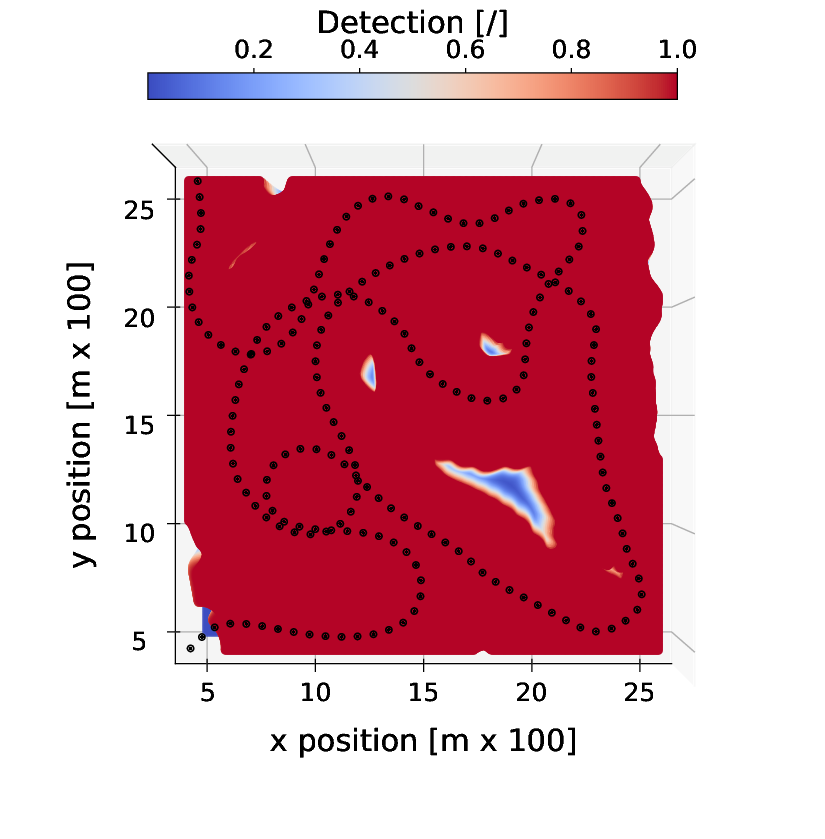}
\end{minipage}
\hspace{.5cm}
\begin{minipage}{0.33\textwidth}
\includegraphics[scale=0.22]{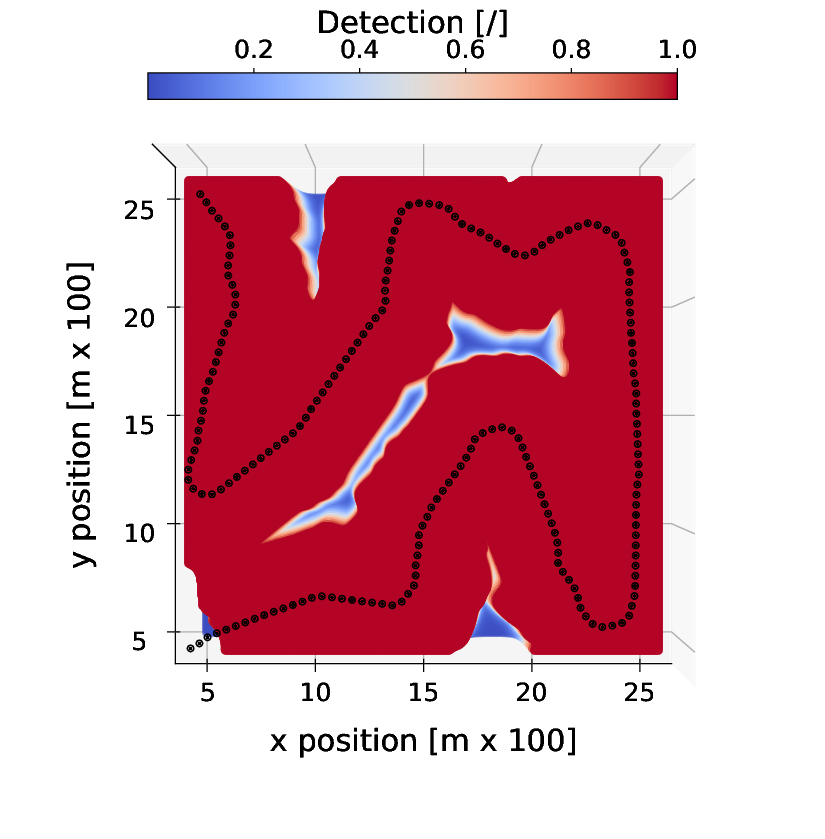}
\end{minipage}
\caption{Three solutions superimposed on the integration points and corresponding detection: (Left) a solution obtained while using a Rank-1 Lattice rule, (Middle) a solution obtained while using a set of MC points, (Right) a solution obtained while using a different set of MC points. The blue line is the trajectory of the autonomous vehicle, and the orange points represent the integration points. The color red represents the investigated area, blue the area that has not been investigated, and the black dotted line is the trajectory of the autonomous vehicle.}
\label{Fig:sol_points_traject}
\end{figure}

We observe that the individual results for the residual risk of the qMC point sets are less spread out than the ones from the MC point set. We hypothesise that this `less spreadoutness' in the risk results comes from the spatial regularity of the distribution of the qMC points, and the Rank-1 Lattice points in particular. When using MC points, negative spaces exists where no points are present. The stochastic optimization software does not consider these negative spaces to be part of the domain, and thus does not plan a trajectory through them. This leads to a longer computation and path times in order to satisfy the requested residual risk. This is illustrated in Fig.\,\ref{Fig:sol_points_traject} and Tab.\,\ref{Tab:res_example_3_sol} where we present the results of three simulations whilst using a Rank-1 Lattice point set and Monte Carlo points, superimposed on the trajectories. This spatial regularity will result in our relaxation strategy requiring fewer successive domain increases to satisfy the requested risk, thereby reducing computational cost.

\begin{table}
\setlength{\tabcolsep}{10pt}
\hspace{-0.5cm}
\scalebox{0.85}{
\begin{tabular}{cccccc}
\toprule
 {} &  {qMC } &   \multirow{2}*{MC}  & \multirow{2}*{Sobol'} & {Sobol'} & {Sobol'}\\
 {} &  {Rank-1 Lattice} & {} &{}& Interlaced factor 2& Interlaced factor 3 \\

 \cmidrule(rl{4pt}){1-6}  
 Median & 4.37\,\% &   4.10\,\% & 4.36\,\%&4.18\,\%&4.29\,\%  \\
 Mean & 4.28\,\% & 3.87\,\%  &4.24\,\%&4.00\,\%& 4.17\,\% \\
 Standard deviation & 0.52\,\%&  0.87\,\%    &0.57\,\%&0.78\,\%& 0.61\,\%    \\
 Minimum &  2.14\,\%  & 0.77\,\%&1.81\,\%&0.94\,\%&1.56\,\%\\
 Maximum &  4.99\,\% & 4.99\,\% &4.99\,\%&4.99\,\%&4.99\,\%\\
\bottomrule
\end{tabular}}
\caption{Statistical results of the residual risk belonging to Fig.\,\ref{fig:residual_mcm_risk_adaptive}.}
\label{Tab:residual_mcm_risk_adaptive}
\end{table}

\begin{table}
\setlength{\tabcolsep}{10pt}
\hspace{-0.5cm}
\scalebox{0.85}{
\begin{tabular}{cccccc}
\toprule
 {} &  {qMC } &   \multirow{2}*{MC}  & \multirow{2}*{Sobol'} & {Sobol'} & {Sobol'}\\
 {} &  {Rank-1 Lattice} & {} &{}& Interlaced factor 2& Interlaced factor 3 \\

 \cmidrule(rl{4pt}){1-6}  
 Median & 139.13\,s &  307.59\,s & 156.55\,s & 237.62\,s & 173.57\,s\\
 Mean & 137.41\,s & 312.52\,s & 156.23\,s & 235.02\,s& 167.73\,s\\
 Standard deviation & 44.88\,s &  138.32\,s   &54.01\,s &89.07\,s &58.18\,s   \\
 Minimum &  13.62\,s & 13.37\,s&14.24\,s&14.77\,s& 14.39\,s\\
 Maximum &  298.74\,s & 1069.41\,s &344.17\,s&526.73\,s&335.62\,s\\
\bottomrule
\end{tabular}}
\caption{Statistical results of the computation time belonging to  Fig.\,\ref{fig:adaptive_computation_time}}
\label{Tab:adaptive_computation_time}
\end{table}

\begin{table}
\setlength{\tabcolsep}{10pt}
\hspace{-0.5cm}
\scalebox{0.85}{
\begin{tabular}{cccccc}
\toprule 
 {} &  {qMC } &   \multirow{2}*{MC}  & \multirow{2}*{Sobol'} & {Sobol'} & {Sobol'}\\
 {} &  {Rank-1 Lattice} & {} &{}& Interlaced factor 2& Interlaced factor 3 \\

 \cmidrule(rl{4pt}){1-6}  
 Median & 2936\,s &  3475\,s & 3011\,s&3263\,s&3063\,s  \\
 Mean & 3225\,s & 3781\,s & 3305\,s & 3598\,s &3355\,s\\
 Standard deviation & 774\,s&  924\,s & 788\,s & 889\,s &798\,s        \\
 Minimum &  2540\,s  & 2668\,s& 2568\,s&2645\,s&2576\,s\\
 Maximum &  8482\,s & 9987\,s&7461\,s&9435\,s&7854\,s \\
\bottomrule
\end{tabular}}
\caption{Statistical results of the path time belonging to  Fig.\,\ref{fig:adaptive_run_time}}
\label{Tab:adaptive_run_time}
\end{table} 

\begin{table}
\centering
\scalebox{1.0}{
\begin{tabular}{cccc}
\toprule
 {} &  {qMC } &   \multirow{2}*{MC - 1} & \multirow{2}*{MC - 2} \\
 {} &  {Rank-1 Lattice} &    \\
 \cmidrule(rl{4pt}){1-4}  
 Maximal Requested Residual MCM Risk & 0.05\,\% &   0.05\,\% &  0.05\,\% \\
 Obtained Residual MCM Risk & 0.043\,\% & 0.014\,\% &  0.043\,\% \\
 Path time & 2677\,s  & 5424\,s & 3431\,s\\
 Initial Domain & $\left[5.0,25.0\right]^2$&  $\left[5.0,25.0\right]^2$    &  $\left[5.0,25.0\right]^2$   \\
 End Domain &  $\left[4.4,25.6\right]^2$  &  $\left[4.4,25.6\right]^2$ & $\left[3.0,27.0\right]^2$\\
\bottomrule
\end{tabular}}
\caption{Most relevant results pertaining to the solutions presented in Fig.\,\ref{Fig:sol_points_traject}.}
\label{Tab:res_example_3_sol}
\end{table}


\subsection{Relaxation strategy for convex quadrilateral domains}
\begin{figure}
\hspace{1.5cm}
\begin{minipage}{0.45\textwidth}
\includegraphics[scale=0.40]{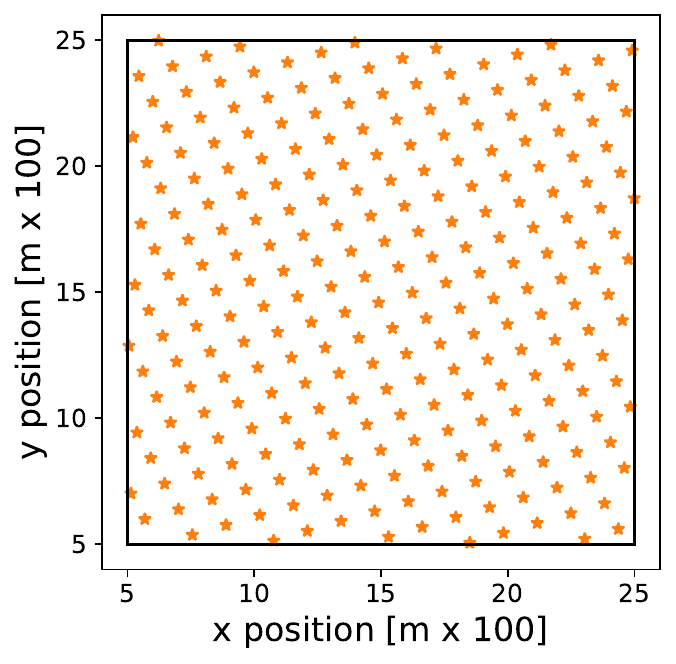}
\end{minipage}
\begin{minipage}{0.50\textwidth}
\includegraphics[scale=0.40]{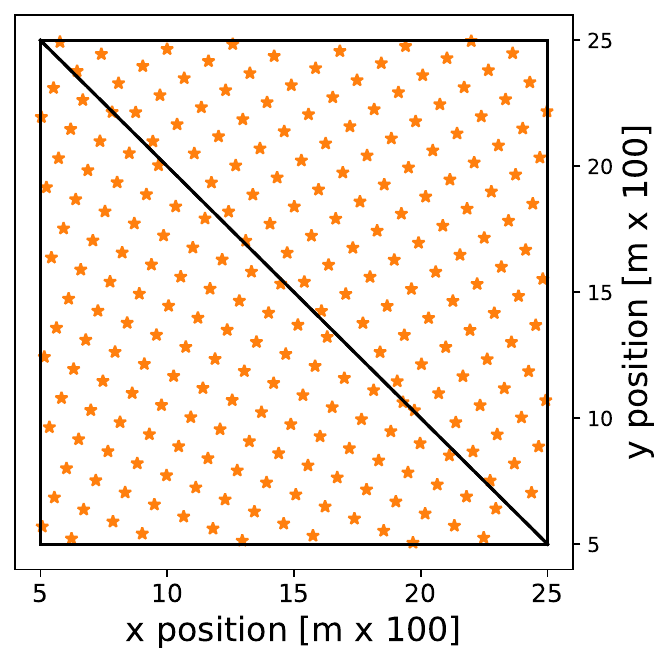}
\end{minipage}
\caption{(Left) qMC points generated on a square domain. (Right) qMC points generated on two triangular domains which make up a square domain.}
\label{fig:triag_vs_quad}
\end{figure}

We now present a comparison in terms of computational cost when combining the relaxation strategy with qMC points generated on a square domain, and with qMC points generated on triangular domains, when considering a square domain $\Omega = \left[5,25\right]^2$. The respective domains we consider and an example of the qMC points are shown in Fig.\,\ref{fig:triag_vs_quad}. We show results for the residual risk, computation time and path time in Fig.\,\ref{fig:1}, Fig.\,\ref{fig:2} and Fig.\,\ref{fig:3}, and in Tab.\,\ref{Tab:1}, Tab.\,\ref{Tab:2}, and Tab.\,\ref{Tab:3}. The same approach  is used as before i.e., we performed 1000 independent simulations with a maximal requested residual risk of $5\%$. For all the considered results, we observe that the performance of quadrilateral qMC points is similar to the performance of triangular qMC points for our MCM implementation. 

\begin{table}
\begin{minipage}[b]{0.4\linewidth}
\hspace{-1.8cm}
\includegraphics[scale=0.25]{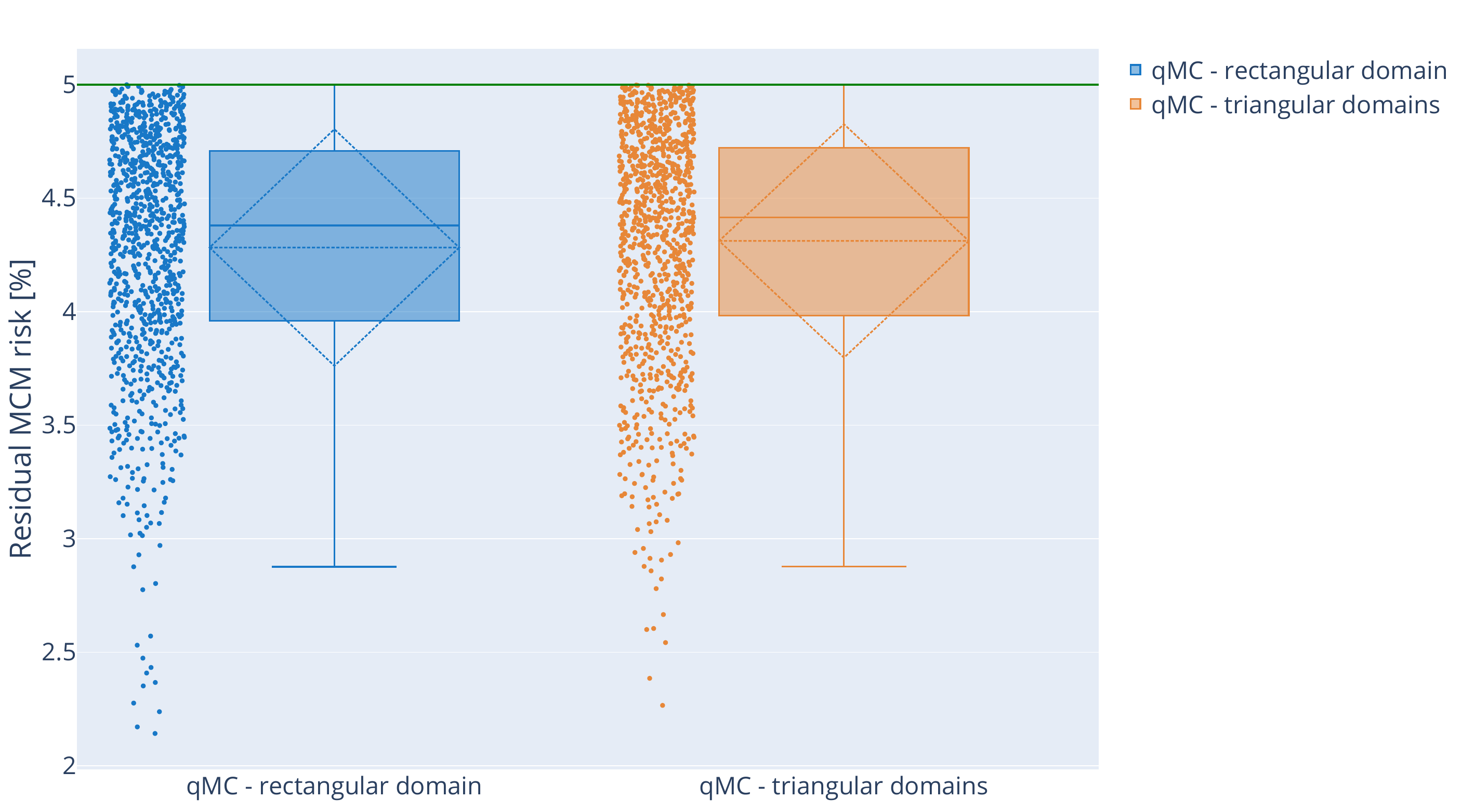}
\captionof{figure}{Comparing the residual MCM risk.}
\label{fig:1}
\end{minipage}
\hfill
\begin{minipage}[b]{0.4\linewidth}
\vspace{-5cm}
\scalebox{0.95}{
\begin{tabular}{ccc}
\toprule
 {} &  {qMC } &   {qMC }  \\
 {} &  {rectangular} &  {triangular}  \\
 \cmidrule(rl{4pt}){1-3}  
 Median & 4.37\,\% &   4.41\,\%  \\
 Mean & 4.28\,\% &   4.31\,\% \\
 Standard deviation & 0.52\,\%&  0.51\,\%        \\
 Minimum &  2.14\,\%  & 2.26\,\%\\
 Maximum &  4.99\,\% & 4.99\,\% \\
\bottomrule
\end{tabular}}
\vspace{1.5cm}
\caption{Numerical results belonging to Fig.\,\ref{fig:1}}
\label{Tab:1}
  \end{minipage}%
\end{table} 

\begin{table}
\begin{minipage}[b]{0.4\linewidth}
\hspace{-1.8cm}
\includegraphics[scale=0.25]{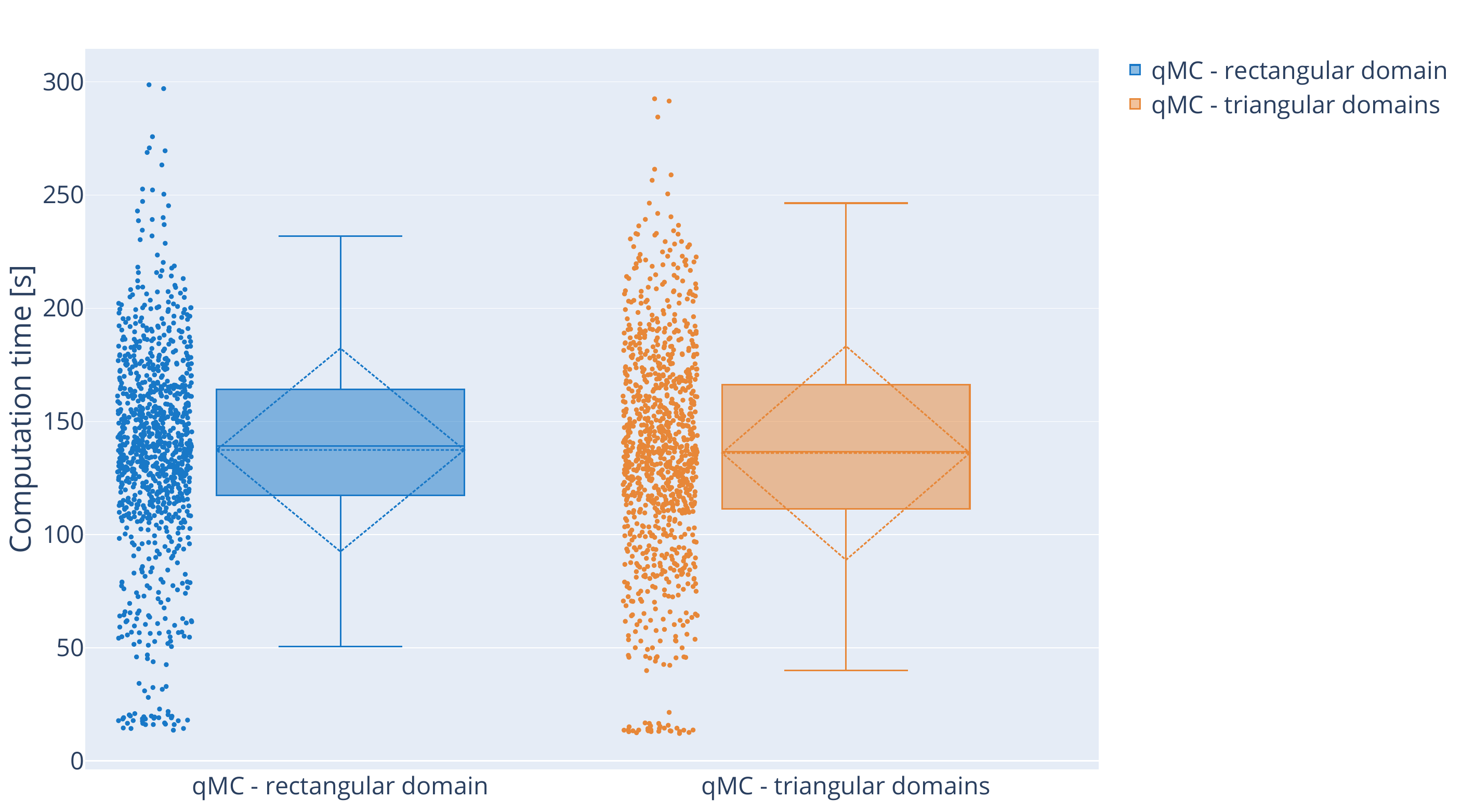}
\captionof{figure}{Comparing the computational times.}
\label{fig:2}
\end{minipage}
\hfill
\begin{minipage}[b]{0.4\linewidth}
\vspace{-5cm}
\scalebox{0.95}{
\begin{tabular}{ccc}
\toprule 
 {} &  {qMC } &   {qMC }  \\
 {} &  {rectangular} &  {triangular}  \\
 \cmidrule(rl{4pt}){1-3}  
 Median & 139.13\,s &   136.56\,s  \\
 Mean & 137.13\,s &   136.06\,s \\
 Standard deviation & 44.88\,s&  47.22\,s        \\
 Minimum &  13.62\,s  & 12.16\,s\\
 Maximum &  298.74\,s & 292.55\,s \\
\bottomrule
\end{tabular}}
\vspace{1.5cm}
\caption{Numerical results belonging to Fig.\,\ref{fig:2}}
\label{Tab:2}
  \end{minipage}%
\end{table} 

\begin{table}
\begin{minipage}[b]{0.4\linewidth}
\hspace{-1.8cm}
\includegraphics[scale=0.25]{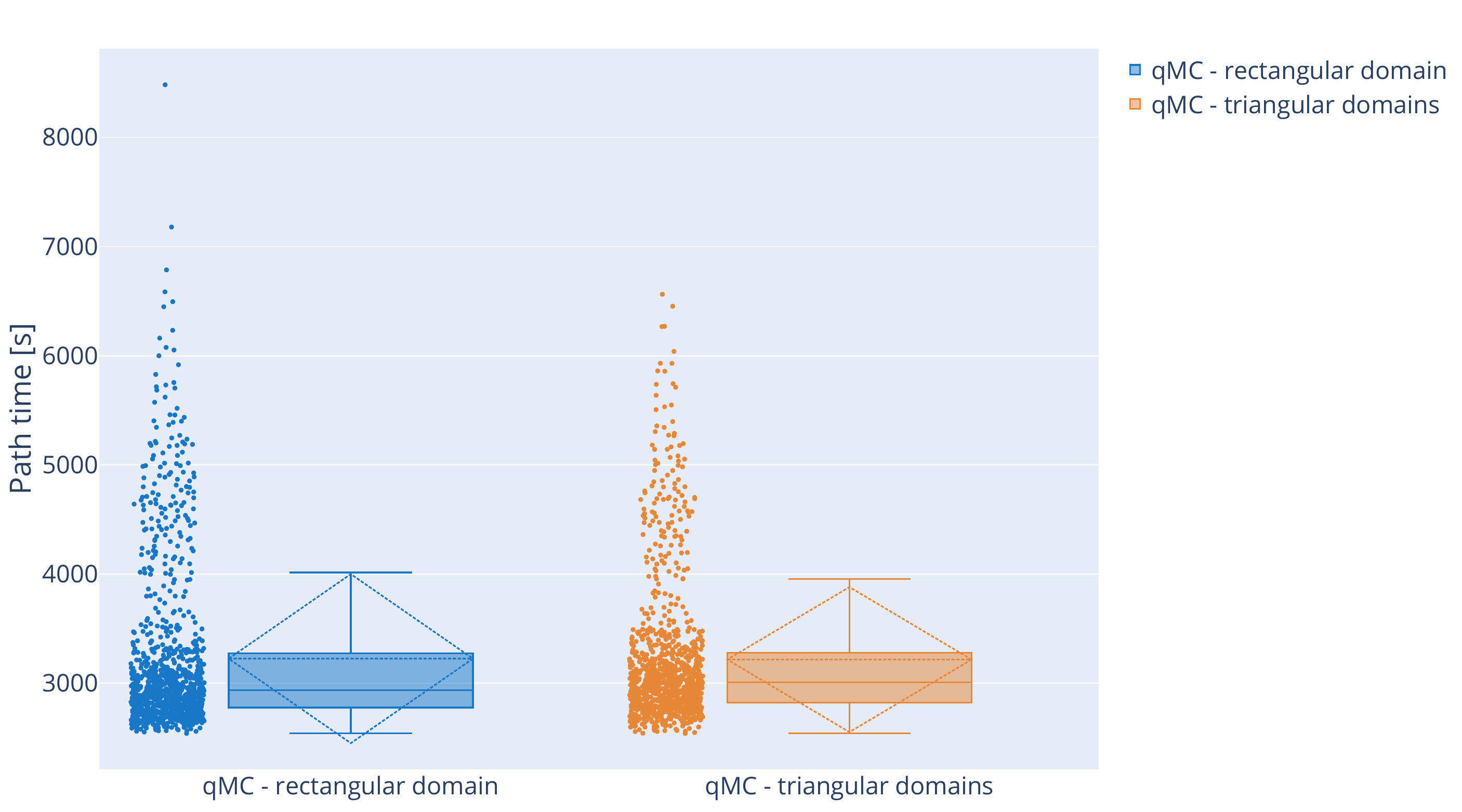}
\captionof{figure}{Comparing the path times.}
\label{fig:3}
\end{minipage}
\hfill
\begin{minipage}[b]{0.4\linewidth}
\vspace{-5cm}
\scalebox{0.95}{
\begin{tabular}{ccc}
\toprule  
 {} &  {qMC } &   {qMC }  \\
 {} &  {rectangular} &  {triangular}  \\
 \cmidrule(rl{4pt}){1-3}  
 Median & 2936\,s &   3007\,s   \\
 Mean & 3225\,s & 3216\,s   \\
 Standard deviation & 774\,s&  664\,s        \\
 Minimum &  2540\,s  & 2539\,s\\
 Maximum &  8482\,s & 6562\,s \\
\bottomrule
\end{tabular}}
\vspace{1.5cm}
\caption{Numerical results belonging to Fig.\,\ref{fig:3}}
\label{Tab:3}
  \end{minipage}%
\end{table} 

\begin{figure}[h]

\begin{minipage}{0.5\textwidth}
\includegraphics[scale=0.5]{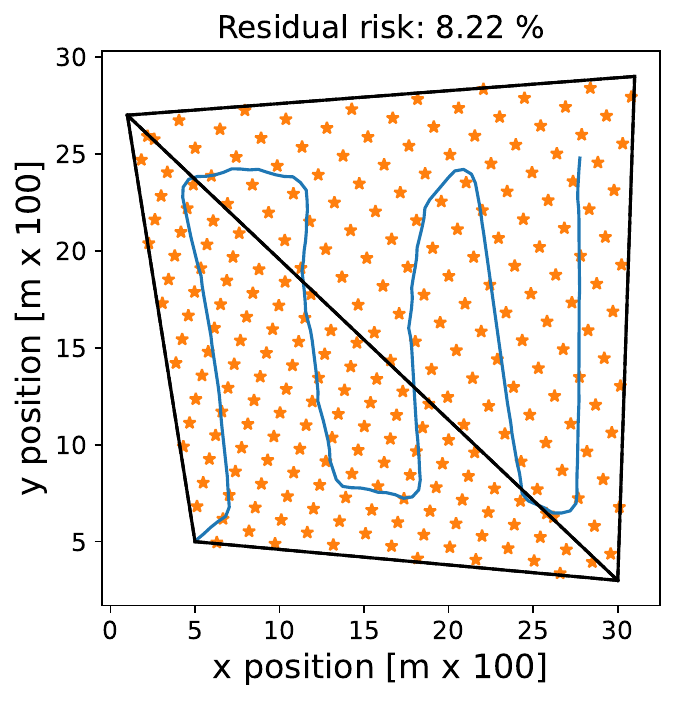}
\end{minipage}
\begin{minipage}{0.5\textwidth}
\includegraphics[scale=0.5]{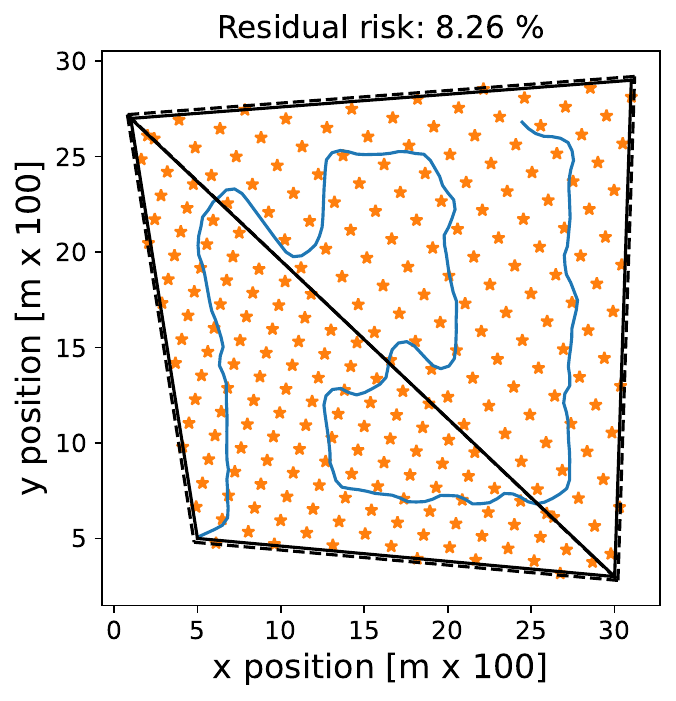}
\end{minipage}

\begin{minipage}{0.5\textwidth}
\includegraphics[scale=0.5]{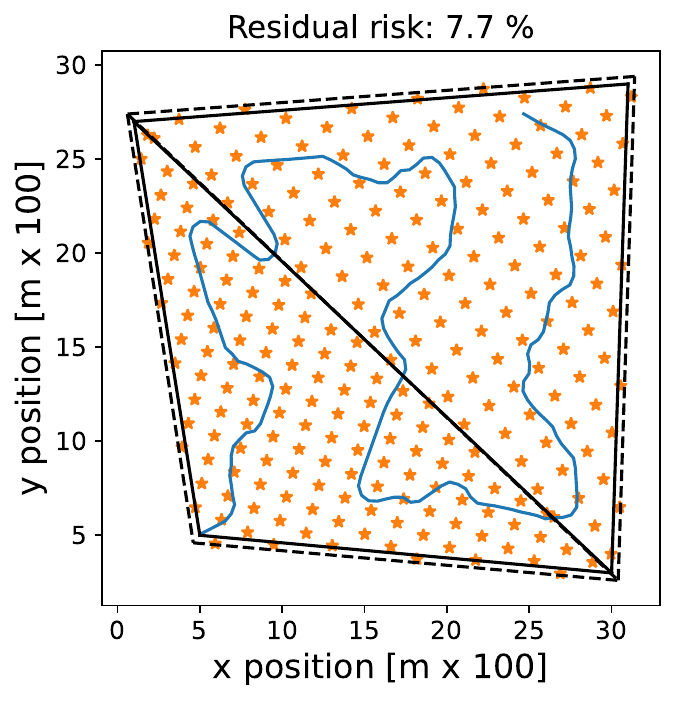}
\end{minipage}
\begin{minipage}{0.5\textwidth}
\includegraphics[scale=0.5]{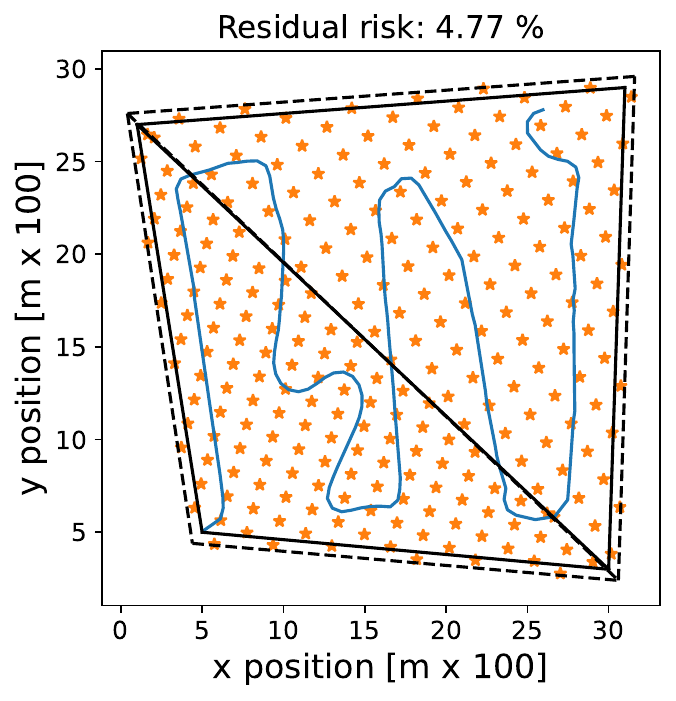}
\end{minipage}

\begin{minipage}{0.99\textwidth}
\hspace{2.5cm}
\includegraphics[scale=0.25]{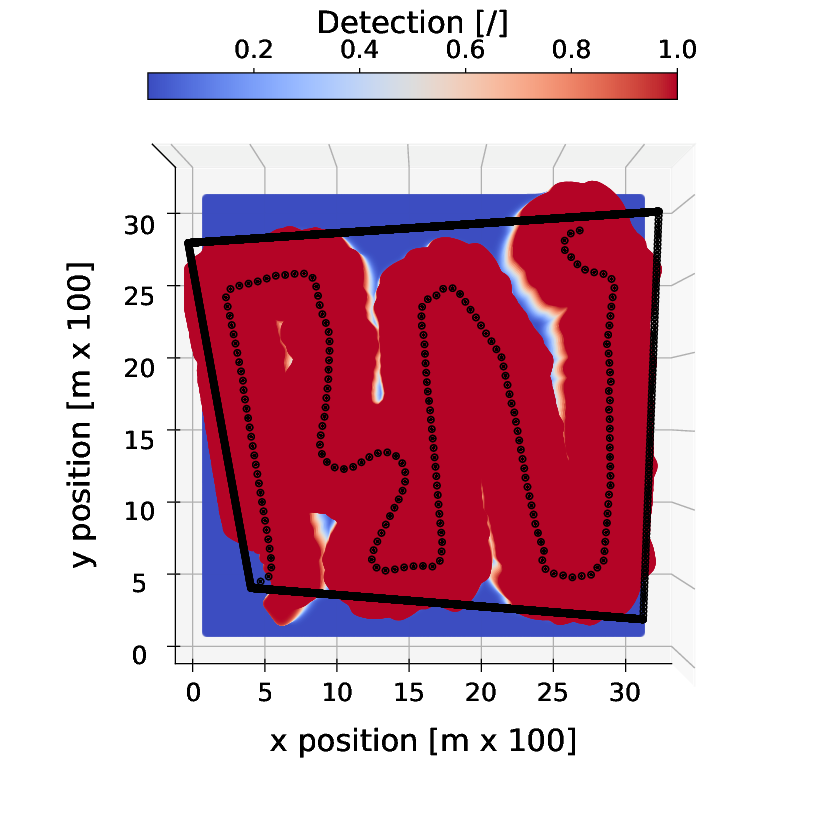}
\end{minipage}
\caption{Illustrative example of our relaxation strategy applied to a convex quadrilateral domain.}
\label{fig:triag_dom_increase}
\end{figure}

Last,  we present a visual representation of our relaxation strategy applied to a convex quadrilateral domain in Fig.\,\ref{fig:triag_dom_increase}. The full black lines represent the initial domain $\Omega_\text{init}$, while the black dotted lines represent the increased domain. We again requested a residual risk of maximally $5\,\%$. Starting with a residual risk of $8.22\,\%$, only three domain increases are necessary to satisfy the requested residual risk.

\section{Conclusion}
In this work we modelled and implemented the mine countermeasure mission problem where the presence of mines is suspected in a user defined domain in a stochastic optimal control framework. One of the most important constraints in our formulation consists of the residual risk integral, which is implemented as an expected value integral. In order to evaluate this integral in the optimization routine we use different qMC and MC point sets. However, we found that at the end of the optimization routine  the requested residual MCM risk is typically not satisfied.  In order to tackle this issue, we proposed and implemented a relaxation strategy. This relaxation consists of successively increasing the square search domain  until the requested residual MCM risk is satisfied. We combined our strategy with different qMC and MC points sets, and found that the Rank-1 Lattice rule yields the best results in terms of computation time. The Rank-1 Lattice rule has a speedup up to factor two when compared to MC. Additionally, we proposed a strategy to simulate convex quadrilateral search domains based on the generation of triangular qMC points. We combined this with our relaxation strategy and showed good results. For future work, we will combine our relaxation strategy with the modelisation of certain oceanographic features such as sand ripples present on the bottom of the ocean. These sand ripples impact the detection of the mines. We also plan to simulate multiple autonomous vehicles in the same search domain, and expand our formulation of convex quadrilateral search domains to domains of more than four sides.
\clearpage

\bibliographystyle{spmpsci}
\bibliography{refs}

\end{document}